\documentclass[12pt,english]{smfart}

\usepackage[francais,english]{babel}

\usepackage{bull, smfthm}

\def\proofname{Proof}
{\hfill\qed\end{trivlist}}
\newenvironment{Proof*}{\begin{trivlist}\item[\hskip%
\labelsep{\bf\proofname.\quad}]}%
{\end{trivlist}}

\def\nextref#1#2#3#4{\bibitem{#1} {\sc {#2}.---} {\it {#3}.} {#4.}}

\def\poi#1#2#3#4#5#6#7{\def\un{#5#6#7}\def\deux{#6#7}
\def\trois{#2#4} \def\cinq{#3#4#5}
\ifx\un\empty {#1}_{#2}{#3}{#1}_{#4} \else
\ifx\deux\empty {#5}(#1_{#2}){#3}{#5}(#1_{#4}) \else
\ifx\trois\empty {#5}_{#6}(#1){#3}{#5}_{#7}(#1) \else
{#5_{#6}}(#1_{#2}){#3}{#5_{#7}}(#1_{#4}) \fi \fi \fi}

\def\dv#1#2{\langle {#1},{#2}\rangle}
    
\def\tk#1#2{{#2}\otimes _{#1}}

\def\ec#1#2#3#4#5{\def\un{#3#4#5}\def\deux{#3#5}\def\trois{#3}
\def\four{#2#4#5}\def\five{#2#5}\def\six{#2}\def\seven{#3#4}
\def\eight{#2#4} \def\nine{#2#3#4}
\ifx\nine\empty {\rm #1}_{#5} \else
\ifx\un\empty {\rm #1}({\goth #2}) \else
\ifx\deux\empty {\rm #1}({\goth #2}_{#4}) \else
\ifx\trois\empty {\rm #1}_{#5}({\goth #2}_{#4}) \else
\ifx\four\empty {\rm #1}(#3) \else
\ifx\five\empty {\rm #1}(#3_{#4}) \else
\ifx\six\empty {\rm #1}_{#5}(#3_{#4}) \else
\ifx\seven\empty {\rm #1}_{#5} ({\goth#2})\else                                               
\ifx\eight\empty {\rm #1}_{#5}({#3})                                               
\fi \fi \fi \fi \fi \fi \fi \fi \fi}
\def\hec#1#2#3#4#5{\def\un{#3#4#5}\def\deux{#3#5}\def\trois{#3}
\def\four{#2#4#5}\def\five{#2#5}\def\six{#2}\def\seven{#3#4}
\def\eight{#2#4} \def\nine{#2#3#4}
\ifx\nine\empty \hat{{\rm #1}}_{#5} \else
\ifx\un\empty \hat{{\rm #1}}({\goth #2}) \else
\ifx\deux\empty \hat{{\rm #1}}({\goth #2}_{#4}) \else
\ifx\trois\empty \hat{{\rm #1}}_{#5}({\goth #2}_{#4}) \else
\ifx\four\empty \hat{{\rm #1}}(#3) \else
\ifx\five\empty \hat{{\rm #1}}(#3_{#4}) \else
\ifx\six\empty \hat{{\rm #1}}_{#5}(#3_{#4}) \else 
\ifx\seven\empty \hat{{\rm #1}}_{#5} ({\goth#2})  \else                                             
\ifx\eight\empty \hat{{\rm #1}}_{#5}({#3}) 
\fi \fi \fi \fi \fi \fi \fi \fi \fi}
\def\e#1#2{\ec {#1}#2{}{}{}}
\def\es#1#2{\ec {#1}{}{#2}{}{}}

\def\ad{{\rm ad}\hskip .1em}

\def\Bbb{\mathbb}
\def\goth{\mathfrak}

\def\j#1#2{\def\deux{#2} \ifx\deux\empty {\rm rk}\hskip .125em{{\goth #1}} \else {\rm rk}\hskip .125em{{\goth #1}_{#2}} \fi}
\def\ji#1#2{\def\deux{#2} \ifx\deux\empty {\rm i}_{{\goth #1}} \else {\rm i}_{{\goth #1}({#2})} \fi}
\def\ai#1#2#3{\def\deux{#2#3} \def\trois{#3} \def\quatre{#2} 
\ifx\deux\empty \es S{{\goth #1}}^{{\goth #1}} \else
\ifx\trois\empty \es S{{\goth #1}(#2)}^{{\goth #1}(#2)} \else
\ifx\quatre\empty \es S{{\goth #1}_{#3}}^{{\goth #1}_{#3}} \else
\es S{{\goth #1}_{#3}(#2)}^{{\goth #1}_{#3}(#2)} \fi \fi \fi}
\def\an#1#2{\def\deux{#2} \ifx\deux\empty {\cal O}_{#1} \else {\cal O}_{#1,#2} \fi }
\def\han#1#2{\def\deux{#2} \ifx\deux\empty {\hat{{\cal O}}}_{#1} \else {\hat{{\cal O}}}_{#1,#2} \fi }

\def\r#1#2#3{\relax \def\un{#3} \relax \ifx\un\empty #1_{{\goth #2}} \else #1_{{\goth #2},{\goth #3}} \fi}
\def\u{\rm u}
\def\dim{{\rm dim}\hskip .125em}
\def\dd{{\rm d}}

\def\ad{{\rm ad}\hskip .1em}

\def\det{{\rm det}\hskip .125em}

\def\r{{\rm r}}

\def\sl2{{\goth s}{\goth l}_{2}({\Bbb C})}

\title[Homogeneous affine varieties]{A remark on homogeneous affine varieties and 
related matters.}

\author[Charbonnel]{Jean-Yves Charbonnel}

\address{Universit\'e Paris 7 - CNRS \\
Institut de Math\'ematiques de Jussieu \\
Th\'eorie des groupes \\
Case 7012 \\ 2 Place Jussieu \\
75251 Paris Cedex 05, France }
\email{jyc@math.jussieu.fr}



\begin{document}

\begin{abstract}
In this note we give an example of affine quotient $G/H$ where $G$ is an affine algebraic
group over an algebraically closed field of characteristic $0$ and $H$ is a unipotent 
subgroup not contained in the unipotent radical of $G$. Some remarks about symmetric
algebras of centralizers of nilpotent elements in simple Lie algebras, in particular 
cases, are added.
\end{abstract}

\maketitle

\section{Introduction.} In this note ${\Bbb K}$ is an algebraically closed field of 
characteristic $0$. Let $G$ be an affine algebraic group over ${\Bbb K}$ and $H$ a closed
subgroup in $G$. We respectively denote by $G_{\u}$ and $H_{\u}$ the unipotent radicals
of $G$ and $H$. It is known \cite{Bi}(Proposition 3) that the homogeneous space $G/H$ is 
affine if and only if it is so for the homogeneous space $G/H_{\u}$. Moreover when 
$H_{\u}$ is contained in $G_{\u}$ then $G/H$ is affine \cite{Bi}(Corollary 2). The 
question to know if the converse is true is an old question. But in this note we give an 
example which shows that it is not true in general. For that purpose we consider a 
simple Lie algebra ${\goth g}$ of type $F_{4}$ over ${\Bbb K}$. In ${\goth g}$ there 
exists a nilpotent element $e$ whose centralizer ${\goth g}(e)$ in ${\goth g}$ has 
dimension $16$ and its reductive factors are simple of dimension $3$. By 
\cite{Ch}(Th\'eor\`eme 3.12), for any $x$ in a non empty open subset of the dual of 
${\goth g}(e)$, its coadjoint orbit is closed. In particular it is an affine
variety. We then show that for any $x$ in a non empty open subset in the dual of
${\goth g}(e)$, the stabilizer ${\goth g}(e)(x)$ of $x$ in ${\goth g}(e)$ is not 
contained in the subset of nilpotent elements of the radical of ${\goth g}(e)$ and any
element of ${\goth g}(e)(x)$ is nilpotent. Furthermore there is a natural question. Let 
${\goth g}$ be a semi-simple Lie algebra and $e$ a nilpotent element. Does exist in the
symmetric algebra of the centralizer of $e$ in ${\goth g}$ a semi-invariant which is not 
an invariant for the adjoint action? We give an example where it is not true and we check
that in this case the index of its centralizer is equal to the rank of ${\goth g}$. In 
the following sections, ${\Bbb K}$ is the ground field and we consider Zariski's topology.
For an algebraic variety $X$ and ${\goth g}$ a Lie algebra acting on $X$ we denote by
${\goth g}(x)$ the stabilizer in ${\goth g}$ of the element $x$ in $X$.  
 
\section{Some remarks about invariants.} \label{i}
Let ${\goth g}$ be a finite dimensional Lie algebra over ${\Bbb K}$. We denote by 
${\goth g}^{*}$ its dual, $\e Sg$ the symmetric algebra of ${\goth g}$, $\ai g{}{}$ the
subalgebra of invariant elements in $\e Sg$ for the adjoint action. We consider the 
coadjoint action of ${\goth g}$ in ${\goth g}^{*}$. For $p$ in $\e Sg$, the differential 
$\dd p(x)$ of $p$ at $x$ is a linear form on ${\goth g}^{*}$. So $\dd p$ is a polynomial
map from ${\goth g}^{*}$ to ${\goth g}$, that is to say $\dd p$ is an element in
$\tk {{\Bbb C}}{\e Sg}{\goth g}$.

\begin{lemm}\label{li}
Let $p$ be an element in $\e Sg$ and ${\goth a}$ an ideal in ${\goth g}$.

{\rm i)} The element $p$ is invariant for the adjoint action of ${\goth a}$ in $\e Sg$ 
if and only if $\dd p(x)$ belongs to the stabilizer in ${\goth g}$ of the restriction
of $x$ to ${\goth a}$ for any $x$ in a non empty open subset in ${\goth g}^{*}$.

{\rm ii)} The elements $p$ is in $\e Sa$ if and only if $\dd p(x)$ is in ${\goth a}$
for any $x$ in a non empty open subset in ${\goth g}^{*}$.
\end{lemm}

\begin{proof}
i) For any $v$ in ${\goth a}$ and any $x$ in ${\goth g}^{*}$, $[v,p](x)$ is equal
to $\dv {\dd p(x)}{v.x}$ where the coadjoint action of $v$ on $x$ is denoted by $v.x$.
But by befinition of the coadjoint action, $\dv {\dd p(x)}{v.x}$ is equal to
$\dv {[\dd p(x),v]}x$. Hence $p$ is invariant for the adjoint action of ${\goth a}$ if 
and only if $\dd p(x)$ belongs to the stabilizer in ${\goth g}$ of the restriction of $x$
to ${\goth a}$ for any $x$ in ${\goth g}^{*}$. So the statement is true because any non 
empty open subset in ${\goth g}^{*}$ is everywhere dense in ${\goth g}^{*}$.

ii) If $p$ is in $\e Sa$, $\dd p$ is in $\tk {{\Bbb C}}{\e Sa}{\goth a}$. Hence for any
$x$ in ${\goth g}^{*}$, ${\goth a}$ contains $\dd p(x)$. Reciprocally if $p$ is not in
$\e Sa$, there exists a basis $\poi x1{,\ldots,}{n}{}{}{}$ in ${\goth g}$ such that 
$x_{1}$ is not in ${\goth a}$ and $\frac{\partial p}{\partial x_{1}}$ is not equal to 
zero. So if $x$ is not a zero of $\frac{\partial p}{\partial x_{1}}$ in ${\goth g}^{*}$, 
$\dd p(x)$ is not in ${\goth a}$.
\end{proof}

As the intersection of two non empty open subsets in ${\goth g}^{*}$ is non empty
the following corollary is an easy consequence of lemma \ref{li}.

\begin{coro}\label{ci}
Let ${\goth a}$ be an ideal in ${\goth g}$. Then ${\goth a}$ does not contain 
${\goth g}(x)$ for any $x$ in a non empty open subset in ${\goth g}^{*}$ if and only if 
$\e Sa$ does not contain $\ai g{}{}$.
\end{coro}

By definition the index $\ji g{}$ of ${\goth g}$ is the smallest dimension of the 
stabilizers for the coadjoint action. By \cite{Dix}(Lemme 7), when ${\goth g}$ is 
algebraic, $\ji g{}$ is the transcendence degree over ${\Bbb K}$ of the field of 
invariants for the adjoint action of ${\goth g}$ in the fraction field of $\e Sg$.

\begin{prop}\label{pi}
We suppose that ${\goth g}$ is algebraic. Then the following conditions are equivalent:

\begin{list}{}{}
\item {\rm 1)} the field of invariants for the adjoint action of ${\goth g}$ in the 
fraction field of $\e Sg$ is the fraction field of $\ai g{}{}$,

\item {\rm 2)} the subalgebra $\ai g{}{}$ contains $\ji g{}$ algebraically independent 
elements,

\item {\rm 3)} for any $x$ in a non empty open subset in ${\goth g}^{*}$, ${\goth g}(x)$
is the image of the map $p\mapsto \dd p(x)$. 
\end{list}
\end{prop}
  
\begin{proof}
We denote by $K$ the fraction field of $\e Sg$ and $K^{{\goth g}}$ the field of 
invariants for the adjoint action of ${\goth g}$ in $K$. We will prove the implications 
$(1) \Rightarrow (2)$, $(2) \Rightarrow (3)$, $(3) \Rightarrow (2)$, 
$(2) \Rightarrow (1)$.

$(1) \Rightarrow (2)$   

As the transcendence degree of $K^{{\goth g}}$ over ${\Bbb K}$ is equal to $\ji g{}$, 
$\ai g{}{}$ contains $\ji g{}$ algebraically independent elements because $K^{{\goth g}}$
is the fraction field of $\ai g{}{}$ by hypothesis.

$(2) \Rightarrow (3)$   

Let $\poi p1{,\ldots,}{\ji g{}}{}{}{}$ be algebraically independent elements in 
$\ai g{}{}$. Then for any $x$ in a non empty open subset in ${\goth g}^{*}$ the elements 
$\poi x{}{,\ldots,}{}{\dd p}{1}{\ji g{}}$ are linearly independent. But for any $x$ in a
non empty open subset in ${\goth g}^{*}$, ${\goth g}(x)$ has dimension $\ji g{}$. Hence
by lemma \ref{li}, (i), ${\goth g}(x)$ is the image of $\ai g{}{}$ by the map
$p\mapsto \dd p(x)$.

$(3) \Rightarrow (2)$   

There exist $x$ in ${\goth g}^{*}$ and $\poi p1{,\ldots,}{\ji g{}}{}{}{}$ such that
$\poi x{}{,\ldots,}{}{\dd p}{1}{\ji g{}}$ are linearly independent. Hence 
$\poi p1{,\ldots,}{\ji g{}}{}{}{}$ are algebraically independent.

$(2) \Rightarrow (1)$   

As $\ai g{}{}$ is integrally closed in $\e Sg$, $K^{{\goth g}}$ is the fraction field
of $\ai g{}{}$.
\end{proof}

\section{About semi-invariants.} \label{si}
We consider ${\goth g}$ as in \ref{i} An element $p$ in $\e Sg$ is called semi-invariant
if $[v,p]$ is colinear to $p$ for any $v$ in ${\goth g}$. We denote by $G$ be the adjoint
algebraic group of ${\goth g}$, $G_{\u}$ the unipotent radical of $G$, ${\goth g}_{0}$ 
the intersection of the kernels of the weights of the semi-invariant elements in $\e Sg$.
For $v$ in ${\goth g}^{*}$, $G(v)$ denotes the stabilizer of $v$ in $G$ for the coadjoint
action of $G$ in ${\goth g}^{*}$. 

\begin{lemm}\label{lsi}
For any $x$ in a non empty open subset $V$ in ${\goth g}^{*}$, ${\goth g}_{0}$ is the 
subset of elements $v$ in ${\goth g}$ such that $\ad v$ is in the Lie algebra of 
$G(x)[G,G]G_{\u}$. Moreover, for any $x$ in $V$, ${\goth g}_{0}$ contains the stabilizer
in ${\goth g}$ of the restriction of $x$ to ${\goth g}_{0}$.
\end{lemm}

\begin{proof}
By \cite{Dix}(Th\'eor\`eme 3.3), for any $x$ in a non empty open subset $V$ in 
${\goth g}^{*}$, ${\goth g}_{0}$ is the subset of elements $v$ in ${\goth g}$ such that 
$\ad v$ is in the Lie algebra of $G(x)[G,G]G_{\u}$. Let $x$ be in $V$, $x_{0}$ its
restriction to ${\goth g}_{0}$, $G(x_{0})$ the stabilizer of $x_{0}$ in $G$. If $T$ is
a maximal torus contained in $G(x_{0})$, the affine subspace $x+{\goth g}_{0}^{\perp}$
is stable for the coadjoint action of $T$ in ${\goth g}^{*}$. So this action has a fixed
point. But any point in ${\goth g}_{0}^{\perp}$ is fix for this action. Hence $T$ is 
contained in $G(x)$ and $G_{0}$. As $G/G_{0}$ is a torus, $G(x_{0})$ is contained in
$G_{0}$ and ${\goth g}_{0}$ contains the stabilizer of $x_{0}$ in ${\goth g}$.
\end{proof}

\begin{coro}\label{csi}
The subalgebra $\ai g{}0$ is generated by the semi-invariant elements in $\e Sg$.
\end{coro}
 
\begin{proof}
Let ${\goth A}$ be the subalgebra generated by semi-invariant elements in $\e Sg$. By 
lemmas \ref{li} and \ref{lsi}, any semi-invariant is contained in $\ec Sg{}0{}$ because
it is invariant for the adjoint action of ${\goth g}_{0}$. So $\ai g{}0$ contains 
${\goth A}$. Any element in $\ai g{}0$ is a sum of weight vectors for the action of $G$ 
in $\ec Sg{}0{}$ because any element in $\ai g{}0$ is invariant by $[G,G]$ and $G_{\u}$. 
So ${\goth A}$ is equal to $\ai g{}0$.
\end{proof}

\begin{rema}\label{rsi}
When ${\goth g}$ is algebraic, $\ad {\goth g}$ is the Lie algebra of $G$ and the Lie 
algebra of $G_{\u}$ is the image by the adjoint representation of the subset 
${\goth g}_{\u}$ of elements $v$ in ${\goth g}$ such that $\ad v$ is in the Lie algebra
of $G_{\u}$. So for any $x$ in a non empty open subset in ${\goth g}^{*}$, 
${\goth g}_{0}$ is the sum of ${\goth g}_{\u}$, $[{\goth g},{\goth g}]$, ${\goth g}(x)$.
\end{rema}

We recall that $K$ is the fraction field of $\e Sg$ and $K^{{\goth g}}$ is the field of 
invariants for the adjoint action of ${\goth g}$ in $K$.

\begin{lemm}\label{l2si}
The field $K^{{\goth g}}$ is contained in the fraction field of the subalgebra
$\ai g{}0$.
\end{lemm}

\begin{proof}
Let $p$ be in $K^{{\goth g}}$. As $\e Sg$ is a unique factorization domain, $p$ has a 
unique decomposition
$$ a = a_{1}^{m_{1}}\cdots a_{k}^{m_{k}} $$ 
where $\poi a1{,\ldots,}{p}{}{}{}$ are prime elements in $\e Sg$ and 
$\poi m1{,\ldots,}{k}{}{}{}$ are integers. But $p$ is invariant for the adjoint action of
${\goth g}$. So because of the unicity of the decomposition, $\poi a1{,\ldots,}{p}{}{}{}$
are semi-invariant elements in $\e Sg$. Then by corollary \ref{csi}, 
$\poi a1{,\ldots,}{p}{}{}{}$ are in $\ai g{}0{}$.
\end{proof}\section{Some remarks about centralizers of nilpotent elements in a simple Lie algebra.} 
\label{cn}
Let ${\goth g}$ be a finite dimensional simple Lie algebra, $G$ the adjoint group of
${\goth g}$, $e$ a nilpotent element in ${\goth g}$. We identify ${\goth g}$ and its dual
by the Killing form $\dv ..$. By Jacobson-Morozov there exist elements $h$ and $f$ 
in ${\goth g}$ such that $e,h,f$ is an ${\goth s}{\goth l}_{2}$-triple in ${\goth g}$. 
The intersection ${\goth l}$ of ${\goth g}(e)$ and ${\goth g}(f)$ is a reductive factor 
in ${\goth g}(e)$. Moreover the restriction of $\dv ..$ to 
${\goth g}(e)\times {\goth g}(f)$ is non degenerate. So ${\goth g}(f)$ is identified to 
the dual of ${\goth g}(e)$.

\subsection{} Let ${\goth t}$ be a subspace in a Cartan subalgebra in ${\goth l}$. We 
denote by ${\goth a}$ the centralizer of ${\goth t}$ in ${\goth g}$. Then ${\goth a}$ is 
a reductive Lie algebra in ${\goth g}$ which contains $e$, $h$, $f$.  Moreover the
subspaces ${\goth a}(f)$ and $[{\goth t},{\goth g}(f)]$ are respectively orthogonal to 
$[{\goth t},{\goth g}(e)]$ and ${\goth a}(e)]$. We consider the coadjoint
action $(v,x)\mapsto v.x$ of ${\goth g}(e)$ in ${\goth g}(f)$. Let $G(e)_{0}$ be the 
identity component of the centralizer of $e$ in $G$. We denote by $(g,x)\mapsto g.x$ the 
coadjoint action of $G(e)_{0}$ on ${\goth g}(f)$.

\begin{lemm}\label{lcn1}
Let $\tau $ be the map $(g,x)\mapsto g.x$ from $G(e)_{0}\times {\goth a}(f)$ to
${\goth g}(f)$.

{\rm i)} Let $W_{0}$ be the set of elements $x$ in ${\goth a}(f)$ such that the linear 
map $v\mapsto v.x$ from $[{\goth t},{\goth g}(e)]$ to ${\goth g}(f)$ is injective. Then
$W_{0}$ is an open subset in ${\goth a}(f)$. Moreover, for any $x$ in $W_{0}$, 
the map $v\mapsto v.x$ is a linear isomorphism from $[{\goth t},{\goth g}(e)]$
to $[{\goth t},{\goth g}(f)]$ and ${\goth g}(e)(x)$ is equal to ${\goth a}(e)(x)$.

{\rm ii)} For any $x$ in $W_{0}$, the map $\tau $ is a submersion at 
$({\rm id}_{{\goth g}},x)$.

{\rm iii)} The restriction of $\tau $ to $G(e)_{0}\times W_{0}$ is a smooth morphism from
\sloppy \hbox{$G(e)_{0}\times W_{0}$} onto an open subset in ${\goth g}(f)$.

{\rm iv)} When ${\goth t}$ is a Cartan subalgebra in ${\goth l}$, for any simple factor 
${\goth a}_{1}$ in ${\goth a}$, the component of $e$ on ${\goth a}_{1}$ is a 
distinguished nilpotent element in ${\goth a}_{1}$.
\end{lemm}

\begin{proof}
i) Let $x$ be in ${\goth a}(f)$. For any $v$ in ${\goth g}(e)$ and $t$ in ${\goth t}$, 
we have
$$\dv w{[t,v].x} = \dv {[w,[t,v]]}x = \dv {[t,[w,v]]}x = -\dv {[w,v]}{[t,x]} =0$$
for any $w$ in ${\goth a}(e)$. Hence the image of $[{\goth t},{\goth g}(e)]$ by the
map $v\mapsto v.x$ is contained in $[{\goth t},{\goth g}(f)]$. As 
$[{\goth t},{\goth g}(e)]$ and $[{\goth t},{\goth g}(f)]$ have the same dimension,
$W_{0}$ is the set of points $x$ in ${\goth a}(f)$ where the map $v\mapsto v.x$ from
$[{\goth t},{\goth g}(e)]$ to ${\goth g}(f)$ has maximal rank. So $W_{0}$ is an open
subset in ${\goth a}(f)$. As ${\goth t}$ is contained in ${\goth g}(e)(x)$ for any $x$ in 
${\goth a}(f)$, ${\goth g}(e)(x)$ is the direct sum of ${\goth a}(e)(x)$ and its 
intersection with $[{\goth t},{\goth g}(e)]$. Then for any $x$ in $W_{0}$, 
${\goth g}(e)(x)$ is equal to ${\goth a}(e)(x)$.

ii) Let $x$ be in $W_{0}$. As the linear map tangent to $\tau $ at 
$({\rm id}_{{\goth g}},x)$ is the map
$$ \ad {\goth g}(e) \times {\goth a}(f) \rightarrow {\goth g}(f) \mbox{ , } 
(v,y) \mapsto v.x + y \mbox{ ,}$$
$\tau $ is a submersion at $({\rm id}_{{\goth g}},x)$ by (i).

iii) For any $x$ in ${\goth a}(f)$ the subset of elements $g$ in $G(e)_{0}$ such that 
$\tau $ is a submersion at $(g,x)$ is stable by left multiplication. So the restriction 
of $\tau $ to $G(e)_{0}\times W_{0}$ is a smooth morphism from $G(e)_{0}\times W_{0}$
onto an open subset in ${\goth g}(f)$.

iv) Let us suppose that ${\goth t}$ is a Cartan subalgebra in ${\goth l}$. Let 
${\goth a}_{1}$ be a simple factor in ${\goth a}$ and $e_{1}$ the component of
$e$ on ${\goth a}_{1}$. As ${\goth t}$ is a Cartan subalgebra in ${\goth l}$, 
${\goth a}_{1}(e_{1})$ has no semi-simple elements because ${\goth a}_{1}(e_{1})$ is the
intersection of ${\goth g}(e)$ and ${\goth a}_{1}$. So $e_{1}$ is a distinguished 
nilpotent element. 
\end{proof}

\begin{coro}\label{ccn1}
The open subset $W_{0}$ is not empty if and only if for any $x$ in a non empty open subset
in ${\goth g}(f)$ the orbit $G(e)_{0}.x$ contains an element whose stabilizer contains
${\goth t}$. Moreover in this case ${\goth g}(e)$ and ${\goth a}(e)$ have the same index.
\end{coro}

\begin{proof}
If $W_{0}$ is not empty then by lemma \ref{lcn1}, (iii), for any $x$ in a non empty open
subset in ${\goth g}(f)$, ${\goth g}(e)(x)$ contains a subalgebra conjugated to 
${\goth t}$ by adjoint action because ${\goth g}(e)(y)$ contains ${\goth t}$ for any $y$ 
in $W_{0}$. Reciprocally we suppose that there exists a non empty open subset $U$ in 
${\goth g}(f)$ such that $G(e)_{0}.x$ contains an element whose stabilizer contains 
${\goth t}$. As we can suppose $U$ invariant by $G(e)_{0}$ there exists an element 
$x$ in ${\goth a}(f)$ which is regular as a linear form on ${\goth g}(e)$ and 
${\goth a}(e)$. In particular ${\goth a}(e)(x)$ and ${\goth g}(e)(x)$ are commutative 
subalgebra. Hence ${\goth g}(e)(x)$ is contained in ${\goth a}(e)(x)$ because it contains
${\goth t}$. So the map $v\mapsto v.x$ from $[{\goth t},{\goth g}(e)]$ to ${\goth g}(f)$ 
is injective and $W_{0}$ contains $x$. Moreover for any $v$ in ${\goth a}(e)(x)$, $x$ is 
orthogonal to $[v,[{\goth t},{\goth g}(e)]]$ because it is contained in 
$[{\goth t},{\goth g}(e)]$. So ${\goth g}(e)(x)$ is equal to ${\goth a}(e)(x)$ and
${\goth a}(e)$ has the same index as ${\goth g}(e)$. 
\end{proof}

\subsection{} For any $\lambda $ in ${\goth t}^{*}$, we respectively denote by 
$E_{\lambda }$ and $F_{\lambda }$ the weight subspaces of weight $\lambda $ for the 
adjoint action of ${\goth t}$ in ${\goth g}(e)$ and ${\goth g}(f)$.

\begin{lemm}\label{lcn2}
Let $\lambda $ be in ${\goth t}^{*}$. Then the subspaces $E_{\lambda }$, $E_{-\lambda }$,
$F_{\lambda }$, $F_{-\lambda }$ have the same dimension. Moreover the eigenvalues of
the restrictions of $\ad h$ to $E_{\lambda }$ and $E_{-\lambda }$ are the same with the 
same multiplicities.
\end{lemm}

\begin{proof}
We denote by ${\goth s}$ the subspace generated by $e,h,f$. It is well known that 
the restriction of $\dv ..$ to ${\goth g}(e)\times {\goth g}(f)$ is non degenerate. But 
for $\mu $ in ${\goth t}^{*}$, $F_{\mu }$ is orthogonal to $E_{\lambda }$ if $\mu $ is 
different from $-\lambda $. So the restriction of $\dv ..$ to 
$E_{\lambda }\times F_{-\lambda }$ is non degenerate because ${\goth g}(e)$ and 
${\goth g}(f)$ are respectively the sum of the subspaces $E_{\mu }$ and $F_{\mu }$ where 
$\mu $ is in ${\goth t}^{*}$. So the dimension of $E_{\lambda }$ and $F_{-\lambda }$ are 
equal. As ${\goth t}$ centralizes $h$, there exists a basis $\poi w1{,\ldots,}{k}{}{}{}$ 
in $F_{-\lambda }$ whose elements are weight vectors for $h$. For $i=1,\ldots,k$, we 
denote by $W_{i}$ the sub-${\goth s}$-module in ${\goth g}$ generated by $w_{i}$. Then 
$W_{i}$ is simple and the intersection of $W_{i}$ and ${\goth g}(e)$ has dimension $1$. 
Let $v_{i}$ be a non zero element in this intersection. As ${\goth t}$ centralizes 
${\goth s}$, any element in $W_{i}$ is a weight vector of weight $-\lambda $ for 
${\goth t}$. So $E_{-\lambda }$ contains $\poi v1{,\ldots,}{k}{}{}{}$. As 
$\poi w1{,\ldots,}{k}{}{}{}$ are linearly independent, the sum of subspaces 
$\poi W1{,\ldots,}{k}{}{}{}$ is direct. So the elements $\poi v1{,\ldots,}{k}{}{}{}$ are 
linearly independent. Then the dimension of 
$E_{-\lambda }$ is bigger than the dimension of $E_{\lambda }$. By the same reasons,
$E_{-\lambda }$ and $F_{\lambda }$ have the same dimension and it is smaller than the 
dimension of $E_{\lambda }$. Hence the subspaces $E_{\lambda }$, $E_{-\lambda }$, 
$F_{\lambda }$, $F_{-\lambda }$ have the same dimension. Moreover 
$\poi v1{,\ldots,}{k}{}{}{}$ is a basis in $E_{-\lambda }$ whose elements are 
eigenvectors for $\ad h$. So if $d$ is an eigenvalue of the restriction of $\ad h$ to
$E_{-\lambda }$ with multiplicity $m$, $-d$ is an eigenvalue of the restriction of 
$\ad h$ to $F_{-\lambda }$ with multiplicity $m$. From the duality between $E_{\lambda }$
and $F_{-\lambda }$ we deduce that $d$ is an eigenvalue of the restriction of $\ad h$ to
$E_{-\lambda }$ with multiplicity $m$.
\end{proof}

Let $\lambda $ be a non zero weight of the adjoint action of ${\goth t}$ in 
${\goth g}(e)$. Let $\poi v1{,\ldots,}{k}{}{}{}$ be a basis in $E_{\lambda }$ and 
$\poi w1{,\ldots,}{k}{}{}{}$ a basis in $E_{-\lambda }$.

\begin{lemm}\label{l2cn2}
Let $M_{\lambda }$ be the matrix
$$M_{\lambda } = \left [ \begin{array}{ccc}
[v_{1},w_{1}] & \cdots & [v_{1},w_{k}] \\
\vdots & \ddots & \vdots \\ { }[v_{k},w_{1}] & \cdots & [v_{k},w_{k}] \end{array} 
\right ] $$
with coefficients in $\es S{{\goth g}(e)}$ and $\delta _{\lambda }$ its determinant.

{\rm i)} The element $\delta _{\lambda }$ is in $\es S{{\goth a}(e)}$. Moreover up to a
multiplicative scalar $\delta _{\lambda }$ does not depend on the choice of the 
basis in $E_{\lambda }$ and $E_{-\lambda }$.

{\rm ii)} For any $x$ in ${\goth a}(f)$, ${\goth g}(e)(x)$ is the sum of its 
intersections with the subspaces $E_{\mu }$ where $\mu $ is in ${\goth t}^{*}$.

{\rm iii)} Let $x$ be in ${\goth a}(f)$. Then the intersection of $E_{\lambda }$ and 
${\goth g}(e)(x)$ is equal to $\{0\}$ if and only if $\delta _{\lambda }(x)$ is different
from $0$.
\end{lemm}

\begin{proof}
By definition, $E_{0}$ and $F_{0}$ are respectively equal to ${\goth a}(e)$ and 
${\goth a}(f)$.

i) For $v$ in $E_{\lambda }$ and $w$ in $E_{-\lambda }$, $[v,w]$ is in ${\goth a}(e)$. So 
$\delta _{\lambda }$ is in $\es S{{\goth a}(e)}$. Let $\alpha $ be a linear automorphism
in $E_{\lambda }$, $\alpha $ its matrix in the basis $\poi v1{,\ldots,}{k}{}{}{}$, 
$\alpha _{i,j}$ the coefficient of $\alpha $ on the $i$-rd line and $j$-rd column. From
the equality:
$$\alpha (v_{i}) = \sum_{j=1}^{k} \alpha _{i,j}v_{j} $$
for $i=1,\ldots,k$, we deduce
$$\alpha M_{\lambda }  = \left [ \begin{array}{ccc}
[\alpha (v_{1}),w_{1}] & \cdots & [\alpha (v_{1}),w_{k}] \\
\vdots & \ddots & \vdots \\ { }[\alpha (v_{k}),w_{1}] & \cdots & [\alpha (v_{k}),w_{k}] 
\end{array} \right ] \mbox{ .}$$
Hence up to a multiplicative scalar $\delta _{\lambda }$ does not depend on the basis
$\poi v1{,\ldots,}{k}{}{}{}$. As it is the same for the basis in $E_{-\lambda }$, the 
proof is done.

ii) Let $x$ be in ${\goth a}(f)$. As ${\goth a}(f)$ is orthogonal to $E_{\mu }$ for 
$\mu $ different from $0$, ${\goth t}$ is contained in ${\goth g}(e)(x)$. Then 
${\goth g}(e)(x)$ is stable by the adjoint action of ${\goth t}$ in ${\goth g}(e)$. So
${\goth g}(e)(x)$ is the sum of its intersections with the subspaces $E_{\mu }$.

iii) Let $v$ be in $E_{\lambda }$ and $\poi a1{,\ldots,}{k}{}{}{}$ its coordinates in 
the basis $\poi v1{,\ldots,}{k}{}{}{}$. The element $v$ is in ${\goth g}(e)(x)$ if and
only if $\dv x{[v,w_{j}]}$ is equal to $0$ for $j=1,\ldots,k$. This condition is 
equivalent to the equality:
$$\left [ \begin{array}{ccc} a_{1} & \cdots & a_{k} \end{array} \right ]
\left [ \begin{array}{ccc} \dv x{[v_{1},w_{1}]} & \cdots & \dv x{[v_{1},w_{k}]} \\
\vdots & \ddots & \vdots \\ \dv x{[v_{k},w_{1}]} & \cdots & \dv x{[v_{k},w_{k}]} 
\end{array} \right ] = 0 \mbox{ .}$$
So the intersection of $E_{\lambda }$ and ${\goth g}(e)(x)$ is equal to $\{0\}$ if and
only if $\delta _{\lambda }(x)$ is different from $0$.
\end{proof}

Let $\delta $ be the product of $\delta _{\lambda }$ where $\lambda $ is a non zero 
weight of the adjoint action of ${\goth t}$ in ${\goth g}(e)$. 

\begin{coro}\label{ccn2}
The element $\delta $ in $\es S{{\goth a}(e)}$ is different from $0$ if and only if for 
any $x$ in a non empty open subset in ${\goth g}(f)$ the orbit $G(e)_{0}.x$ contains an 
element whose stabilizer contains ${\goth t}$.
\end{coro}

\begin{proof}
We consider the open subset $W_{0}$ in ${\goth a}(f)$ introduced in lemma \ref{lcn1}, (i).
Then an element $x$ in ${\goth a}(f)$ is in $W_{0}$ if and only if $\delta (x)$
is different from $0$. Hence the corollary is a consequence of corollary \ref{ccn1}.
\end{proof}

We recall that an element $x$ in ${\goth g}(f)$ is called regular if the dimension of 
${\goth g}(e)(x)$ is minimal. Moreover the subset of regular elements in ${\goth g}(f)$
is open.

\begin{coro}\label{c2cn2}
If $\delta $ is equal to $0$, then for any regular element $x$ in ${\goth g}(f)$, 
${\goth g}(e)(x)$ does not contain a conjugate of ${\goth t}$ by the adjoint group of 
${\goth g}(e)$.
\end{coro}

\begin{proof}
Let us suppose that there exists a regular element $x$ in ${\goth g}(f)$ such that 
${\goth g}(e)(x)$ contains a conjugate of ${\goth t}$ by the adjoint group of 
${\goth g}(e)$. As the subset of regular elements in ${\goth g}(f)$ is stable for
the action of the adjoint group of ${\goth g}(e)$ we can suppose that ${\goth g}(e)(x)$ 
contains ${\goth t}$. As $x$ is regular, ${\goth g}(e)(x)$ is commutative and contained
in ${\goth a}(e)$. Then by lemma \ref{l2cn2}, (iii), $\delta _{\lambda }(x)$ is not
equal to $0$ for any non zero weight $\lambda $ of the adjoint action of ${\goth t}$ in
${\goth g}(e)$.
\end{proof}

We denote by ${\goth g}(e)_{\u}$ the subset of nilpotent elements in the radical of 
${\goth g}(e)$.

\begin{coro}\label{c3cn2}
Let $d$ be the biggest eigenvalue of $\ad h$. We suppose that the kernel of $\ad h -d$
has dimension smaller than $3$ and does not centralize ${\goth t}$. Moreover if 
the kernel of $\ad h-d$ has dimension $3$, ${\goth t}$ is contained in the center of
${\goth l}$.

{\rm i)} For any regular element $x$ in ${\goth g}(f)$, ${\goth g}(e)(x)$ 
does not contain a conjugate of ${\goth t}$ by the adjoint representation.

{\rm ii)} If ${\goth l}$ has rank $1$ then for any $x$ in a non empty open subset in
${\goth g}(f)$, the elements of ${\goth g}(e)(x)$ are nilpotent.

{\rm iii)} If ${\goth t}$ is the center of ${\goth l}$, then the symmetric algebra 
$\es S{{\goth g}(e)}$ of ${\goth g}(e)$ contains semi-invariant elements which are not in
$\ai ge{}$.
\end{coro}

\begin{proof}
i) As $d$ is the biggest eigenvalue of $\ad h$, the kernel of $\ad h - d$ is contained in
${\goth g}(e)$. By hypothesis and lemma \ref{lcn2}, there exists a non zero linear form 
$\lambda $ on ${\goth t}$ such that $\lambda $ and $-\lambda $ are weights of the 
adjoint action of ${\goth t}$ in the kernel $\ad h - d$. Let $w_{1}$ be a non zero 
element in $E_{\lambda }$ such that $[h,w_{1}]$ is equal to $dw_{1}$. As $d$ is the 
biggest eigenvalue of the restriction of $\ad h$ to ${\goth g}(e)_{\u}$, $w_{1}$ 
centralizes ${\goth g}(e)_{\u}$ because any eigenvalue of the restriction of $\ad h$ to 
${\goth g}(e)_{\u}$ is strictly positive. If there exists a non zero element $v$ in the 
intersection of $E_{-\lambda }$ and ${\goth l}$, the dimension of the kernel of 
$\ad h - d$ is equal to $2$ and $[v,w_{1}]$ is an element in the kernel of $\ad h-d$ 
which centralizes ${\goth t}$. So $w_{1}$ centralizes $E_{-\lambda }$ and 
$\delta _{\lambda }$ is equal to $0$. Hence by corollary \ref{c2cn2}, for any regular
element $x$ in ${\goth g}(f)$, ${\goth g}(e)(x)$ does not contain a conjugate of 
${\goth t}$ by the adjoint group of ${\goth g}(e)$. 

ii) We suppose that ${\goth l}$ has rank $1$. Then for any non zero semi-simple element 
$v$ in ${\goth g}(e)$ the line containing $v$ is conjugate to ${\goth t}$ by the adjoint 
action. So by (i), for any $x$ in a non empty open subset the elements of 
${\goth g}(e)(x)$ are nilpotent.

iii) We suppose that ${\goth t}$ is the center of ${\goth l}$. Then by (i), for any $x$ 
in a non empty open subset, ${\goth g}(e)$ strictly contains the sum of the subspaces
${\goth g}(e)(x)$, $[{\goth l},{\goth l}]$, ${\goth g}(e)_{\u}$. Hence by corollary 
\ref{csi}, $\es S{{\goth g}(e)}$ contains semi-invariant elements which are not
invariant.
\end{proof}

\begin{lemm}\label{l3cn2}
Let $d$ be the biggest eigenvalue of $\ad h$. We suppose that the following conditions
are satisfied:
\begin{list}{}{}
\item {\rm 1)} the kernel of $\ad h-d$ has dimension smaller than $3$ and does not 
centralize ${\goth t}$,
\item {\rm 2)} if the kernel of $\ad h - d$ has dimension $3$ then ${\goth t}$ is 
contained in the center of ${\goth l}$,
\item {\rm 3)} if $\mu $ is a non zero weight of the adjoint action of ${\goth t}$ in
${\goth g}(e)$, $\delta _{\mu }$ is not equal to $0$,
\item {\rm 4)} if $\lambda $ is a non zero weight of the adjoint action of ${\goth t}$ in
the kernel of $\ad h-d$, there exists a principal minor in the matrix $M_{\lambda }$ 
which is not equal to $0$. 
\end{list}
Then the index of ${\goth a}(e)$ is bigger than the index of ${\goth g}(e)$.
\end{lemm}

\begin{proof}
For $i$ in $\{-1,+1\}$ we denote by ${\goth k}_{i\lambda }$ the intersection of 
$E_{i\lambda }$ and the kernel of $\ad h -d$. Let ${\goth k}$ be the sum of 
${\goth k}_{\lambda }$ and ${\goth k}_{-\lambda }$. Then by condition (1) and lemma
\ref{lcn1}, ${\goth k}$ has dimension $2$. Moreover by condition (2), ${\goth k}$ is
an ideal in ${\goth g}(e)$. Let ${\goth q}$ be the quotient of ${\goth g}(e)$ by 
${\goth k}$. Then the restriction to ${\goth a}(e)$ of the canonical morphism from 
${\goth g}(e)$ to ${\goth q}$ is injective. So we can identify ${\goth a}(e)$ with its 
image in ${\goth q}$. Then ${\goth a}(e)$ is the centralizer of ${\goth t}$ in 
${\goth q}$. For any weight $\mu $ of the adjoint action of ${\goth t}$ in ${\goth q}$ 
we denote by $E'_{\mu }$ the weight subspace of weight $\mu $. Then $E'_{\mu }$ is the 
image of $E_{\mu }$ by the canonical morphism from ${\goth g}(e)$ to ${\goth q}$. When 
$\mu $ is not equal to $\lambda $ or $-\lambda $, $E_{\mu }$ and $E'_{\mu }$ have the 
same dimension. Otherwise, $E'_{\lambda }$ and $E'_{-\lambda }$ have dimension 
$\dim E_{\lambda }-1$. Moreover by condition (4) if $\poi v1{,\ldots,}{k}{}{}{}$ and
$\poi w1{,\ldots,}{k}{}{}{}$ are basis in $E'_{\lambda }$ and $E'_{-\lambda }$, we have
$$\det \left [ \begin{array}{ccc}
[v_{1},w_{1}] & \cdots & [v_{1},w_{k}] \\
\vdots & \ddots & \vdots \\ { }[v_{k},w_{1}] & \cdots & [v_{k},w_{k}] \end{array} 
\right ] \not = 0 \mbox{ ,}$$
because the intersection of $E_{\lambda }$ and ${\goth k}$ centralizes $E_{-\lambda }$.
Then by corollaries \ref{ccn2} and \ref{ccn1}, ${\goth a}(e)$ and ${\goth q}$ have the 
same index. Let $x$ be an element in ${\goth a}(f)$ which is a linear regular form on
${\goth a}(e)$ and ${\goth q}$. Then ${\goth g}(e)(x)$ is equal to 
${\goth a}(e)(x)+{\goth k}$. In particular the dimension of ${\goth g}(e)(x)$ is
$\ji ae+2$ and ${\goth g}(e)(x)$ is not commutative because it contains ${\goth t}$ and
${\goth k}$. Hence $x$ is not a regular linear form on ${\goth g}(e)$ and $\ji ge$ is 
smaller than $\ji ae$.
\end{proof}

We finish this section by giving an example of a simple Lie algebra and a nilpotent 
element in it for which the symmetric algebra of its centralizer contains a 
semi-invariant which is not invariant.

Let ${\goth g}$ be a simple Lie algebra of type $E_{7}$. By the tables in \cite{Ca}, 
${\goth g}$ contains a nilpotent element $e$ such that ${\goth g}(e)$ has dimension $35$.
Considering $e,h,f$ as above, ${\goth l}$ is a direct product of a torus ${\goth t}$ of 
dimension $1$ by a simple Lie algebra of dimension $3$. The biggest eigenvalue of $\ad h$
is $6$ and its multiplicity is $3$. Furthermore the kernel of $\ad h-6$ does not 
centralize ${\goth t}$. Then by corollary \ref{c3cn2}, (iii), the symmetric algebra of 
${\goth g}(e)$ contains semi-invariant elements which are not invariant. The centralizer 
${\goth a}$ of ${\goth t}$ has dimension $27$, ${\goth t}$ is the center of ${\goth a}$, 
${\goth a}(e)$ has dimension $9$. Then $[{\goth a},{\goth a}]$ is isomorphic to the 
direct product of ${\goth s}{\goth l}_{4}$, ${\goth s}{\goth l}_{3}$, 
${\goth s}{\goth l}_{2}$. So the index of ${\goth a}(e)$ is $7$. There exists an 
element $t$ in $\ad {\goth t}$ whose eigenvalues are integers. The strictly positive 
values of these integers are $1,2,3,4$ and their respective multiplicities are $6,4,2,1$. 
Moreover the eigenvalues of the restriction of $t$ to the kernel of $\ad h -6$ are
$-2,0,2$. As ${\goth t}$ has dimension $1$, we get elements 
$\delta _{1},\delta _{2},\delta _{3},\delta _{4}$ in $\es S{{\goth a}(e)}$. For $i$ not
equal to $2$, $\delta _{i}$ is different from $0$. Moreover in the matrix whose 
determinant is $\delta _{2}$ there is a principal minor which is not equal to $0$. Hence 
by lemma \ref{l2cn2}, the index of ${\goth g}(e)$ is smaller than $7$. So by  
Vinberg's result, the index of ${\goth g}(e)$ is $7$.\section{An example of affine quotient.} \label{aq}
Let ${\goth g}$ be a simple Lie algebra of type $F_{4}$ and $G$ its adjoint group. 
Following the tables in \cite{Ca}, there exists a nilpotent element $e$ in ${\goth g}$ 
such that ${\goth g}(e)$ has dimension $16$ and the reductive factors of ${\goth g}(e)$ 
are simple of dimension $3$. Then we have the proposition:

\begin{prop}\label{paq}
Let ${\goth g}(e)^{*}$ be the dual of ${\goth g}(e)$. For any $x$ in a non empty open 
subset in ${\goth g}(e)^{*}$, the coadjoint orbit of $x$ is closed in ${\goth g}(e)^{*}$,
the elements of the stabilizer ${\goth g}(e)(x)$ of $x$ in ${\goth g}(e)$ are nilpotent 
and ${\goth g}(e)(x)$ is not contained in the subset of nilpotent elements in the radical
of ${\goth g}(e)$.
\end{prop}

Denoting by $G^{e}$ the adjoint group of ${\goth g}(e)$, for any $x$ in a non empty
open subset in ${\goth g}(e)^{*}$, the coadjoint orbit $G^{e}.x$ is an affine variety and
the identity component $H$ of the stabilizer of $x$ in $G^{e}$ is a unipotent subgroup,
not contained in the unipotent radical of $G^{e}$. So we get an example of affine 
quotient $G^{e}/H$ where $H$ is a unipotent subgroup in $G^{e}$ not contained in the 
unipotent radical of $G^{e}$.

As above we consider elements $h$ and $f$ in ${\goth g}$ such that $e,h,f$ is an 
${\goth s}{\goth l}_{2}$-triple. We identify ${\goth g}(f)$ and ${\goth g}(e)^{*}$ by the
Killing form. The intersection ${\goth l}$ of ${\goth g}(e)$ and ${\goth g}(f)$ is a 
simple Lie algebra of dimension $3$ and the subset ${\goth g}(e)_{\u}$ of nilpotent 
elements in the radical of ${\goth g}(e)$ is an ideal of dimension $13$. As ${\goth g}(e)$
is a semi-direct product of ${\goth l}$ and ${\goth g}(e)_{\u}$, ${\goth g}(e)$ is a
unimodular Lie algebra. So by \cite{Ch}(Th\'eor\`eme 3.12), for any $x$ in a non empty 
open subset in ${\goth g}(f)$ the coadjoint orbit of $x$ is closed in ${\goth g}(f)$. 
Using computation by logicial Gap4 we find a basis $\poi x1{,\ldots,}{16}{}{}{}$ in 
${\goth g}(e)$ whose image by $\ad h$ is the following sequence
$$ x_{1},0,x_{3},2x_{4},2x_{5},3x_{6},x_{7},4x_{8},3x_{9},4x_{10},5x_{11},4x_{12},
5x_{13},x_{14},0,0 \mbox{ .}$$ 
Moreover $e$ is equal to $x_{4}+x_{5}$. Then $x_{2},x_{15},x_{16}$ is a basis in 
${\goth l}$ and $x_{1},\poi x3{,\ldots,}{14}{}{}{}$ is a basis in ${\goth g}(e)_{\u}$. The
biggest eigenvalue of $\ad h$ is equal to five and its multiplicity is equal to $2$. So
by corollary \ref{c2cn2}, for any $x$ in a non empty open subset in ${\goth g}(f)$,
the elements of ${\goth g}(e)(x)$ are nilpotent. The element $p$ in $\es S{{\goth g}(e)}$
$$ p = -9x_{{13}}x_{{12}}x_{{14}}+9x_{{11}}x_{{8}}x_{{7}}+3/2x_{{13}}x_
{{4}}x_{{6}}-3{x_{{11}}}^{2}x_{{2}}+3{x_{{13}}}^{2}x_{{15}}+3/4x
_{{4}}{x_{{10}}}^{2} $$ $$ \mbox{ }-9/4x_{{5}}{x_{{10}}}^{2}-3x_{{1}}x_{{13}}x_{{
10}}+3x_{{1}}x_{{12}}x_{{11}}-3x_{{3}}x_{{13}}x_{{8}}+3x_{{3}}x_
{{11}}x_{{10}}-3x_{{4}}x_{{12}}x_{{8}} $$ $$ \mbox{ } +9x_{{5}}x_{{12}}x_{{8}}+3/2
x_{{11}}x_{{4}}x_{{9}}+3x_{{11}}x_{{13}}x_{{16}} \mbox{ ,}$$
is invariant for the adjoint action. We remark that $p$ is not contained in
$\es S{{\goth g}(e)_{\u}}$. So by corollary \ref{ci}, for any $x$ in a non empty open
subset in ${\goth g}(f)$, ${\goth g}(e)(x)$ is not contained in ${\goth g}(e)_{\u}$. As
the intersection of a finitely many non empty open subsets in ${\goth g}(f)$ is non 
empty, for any $x$ in a non empty open subset in ${\goth g}(f)$ the adjoint orbit of
$x$ is closed in ${\goth g}(f)$, ${\goth g}(e)(x)$ is not contained in ${\goth g}(e)_{\u}$
and its elements are nilpotent. In appendix a Maple program is given to compute $4$ 
algebraic independent elements $p(1),p(2),p(3),p(4)$ in $\ai ge{}$. The element $p$ above
is equal to $p(3)$. As a consequence of lemma \ref{l3cn2}, the index of ${\goth g}(e)$
is equal to $4$. As there are sufficiently many variables with degree $1$ in the 
polynomials $p(1),p(2),p(3),p(4)$, it is easy to see that the generic fiber of the 
morphism $\tau $ whose comorphism is the canonical injection from
${\Bbb C}[p(1),p(2),p(3),p(4)]$ to $\es S{{\goth g}(e)}$, is irreducible. Hence the 
subfield generated by $p(1),p(2),p(3),p(4)$ is the subfield of invariants for the adjoint
action of ${\goth g}$ in the fraction field of $\es S{{\goth g}(e)}$. The main point 
is that the nullvariety of $p(1),p(2),p(3),p(4)$ in ${\goth g}(e)^{*}$ has codimension 
$4$. Then we deduce that $\tau $ is open and surjective. Hence for any $q$ in $\ai ge{}$
the morphism whose comorphism is the canonical injection from 
${\Bbb C}[p(1),p(2),p(3),p(4)]$ to ${\Bbb C}[p(1),p(2),p(3),p(4)][q]$ is surjective and
quasi finite. So by main Zariski's theorem, this morphism is an isomorphism. Hence 
$\ai ge{}$ is a polynomial algebra generated by $p(1),p(2),p(3),p(4)$. Moreover 
$\es S{{\goth g}(e)}$ is a faithfully flat extension of $\ai ge{}$. The same method proves
that for any distinguished nilpotent element $x$ of an exceptional simple Lie algebra
${\goth e}$, the algebra $\ai ex{}$ is a polynomial algebra and $\es S{{\goth e}(x)}$ is
a faithfully flat extension of $\ai ex{}$.\appendix 
\section{Appendix.} 

In this appendix we give the Maple program which computes four algebraically independent
elements in $\ai g{}e$ when ${\goth g}$ is simple of type F$_{4}$ and $e$ is a nilpotent
element such that ${\goth g}(e)$ has dimension $16$ and its reductive factors are
simple of dimension $3$. In this program $\poi x1{,\ldots,}{d}{}{}{}$ is a basis of
${\goth g}(e)$ and $[x_{i},x_{j}]$ is equal to $L(i,j)$. Moreover for $i=1,\ldots,16$,
$y(i)$ is equal to $[h,x_{i}]$ and the center ${\goth z}$ of ${\goth g}(e)$ has dimension
$1$.

\begin{verbatim}

# The procedures k and f are elements in the kernel of 
  the matrix A.

 with(linalg):

  L := proc(i,j)
    if i=1 and j=2 then return 2*x[3]:
  elif i=1 and j=3 then return 3*x[5]:
  elif i=1 and j=4 then return 0:
  elif i=1 and j=5 then return 0:
  elif i=1 and j=6 then return 2*x[8]:
  elif i=1 and j=7 then return 0:
  elif i=1 and j=8 then return 0:
  elif i=1 and j=9 then return -2*x[10]:
  elif i=1 and j=10 then return 2*x[11]:
  elif i=1 and j=11 then return 0:
  elif i=1 and j=12 then return x[13]:
  elif i=1 and j=13 then return 0:
  elif i=1 and j=14 then return 0:
  elif i=1 and j=15 then return 3*x[14]:
  elif i=1 and j=16 then return x[1]:
  elif i=2 and j=3 then return -3*x[7]:
  elif i=2 and j=4 then return 0:
  elif i=2 and j=5 then return 0:
  elif i=2 and j=6 then return x[9]:
  elif i=2 and j=7 then return 0:
  elif i=2 and j=8 then return x[10]:
  elif i=2 and j=9 then return 0:
  elif i=2 and j=10 then return 2*x[12]:
  elif i=2 and j=11 then return -x[13]:
  elif i=2 and j=12 then return 0:
  elif i=2 and j=13 then return 0:
  elif i=2 and j=14 then return -x[1]:
  elif i=2 and j=15 then return x[16]:
  elif i=2 and j=16 then return -2*x[2]:
  elif i=3 and j=4 then return 0:
  elif i=3 and j=5 then return 0:
  elif i=3 and j=6 then return -2*x[10]:
  elif i=3 and j=7 then return 0:
  elif i=3 and j=8 then return x[11]:
  elif i=3 and j=9 then return 2*x[12]:
  elif i=3 and j=10 then return 2*x[13]:
  elif i=3 and j=11 then return 0:
  elif i=3 and j=12 then return 0:
  elif i=3 and j=13 then return 0:
  elif i=3 and j=14 then return 0:
  elif i=3 and j=15 then return 2*x[1]:
  elif i=3 and j=16 then return -x[3]:
  elif i=4 and j=5 then return 0:
  elif i=4 and j=6 then return 2*x[11]:
  elif i=4 and j=7 then return 0:
  elif i=4 and j=8 then return 0:
  elif i=4 and j=9 then return -2*x[13]:
  elif i=4 and j=10 then return 0:
  elif i=4 and j=11 then return 0:
  elif i=4 and j=12 then return 0:
  elif i=4 and j=13 then return 0:
  elif i=4 and j=14 then return 0:
  elif i=4 and j=15 then return 0:
  elif i=4 and j=16 then return 0:
  elif i=5 and j=6 then return -2*x[11]:
  elif i=5 and j=7 then return 0:
  elif i=5 and j=8 then return 0:
  elif i=5 and j=9 then return 2*x[13]:
  elif i=5 and j=10 then return 0:
  elif i=5 and j=11 then return 0:
  elif i=5 and j=12 then return 0:
  elif i=5 and j=13 then return 0:
  elif i=5 and j=14 then return 0:
  elif i=5 and j=15 then return 0:
  elif i=5 and j=16 then return 0:
  elif i=6 and j=7 then return -2*x[12]:
  elif i=6 and j=8 then return 0:
  elif i=6 and j=9 then return 0:
  elif i=6 and j=10 then return 0:
  elif i=6 and j=11 then return 0:
  elif i=6 and j=12 then return 0:
  elif i=6 and j=13 then return 0:
  elif i=6 and j=14 then return 0:
  elif i=6 and j=15 then return 0:
  elif i=6 and j=16 then return x[6]:
  elif i=7 and j=8 then return x[13]:
  elif i=7 and j=9 then return 0:
  elif i=7 and j=10 then return 0:
  elif i=7 and j=11 then return 0:
  elif i=7 and j=12 then return 0:
  elif i=7 and j=13 then return 0:
  elif i=7 and j=14 then return x[5]:
  elif i=7 and j=15 then return x[3]:
  elif i=7 and j=16 then return -3*x[7]:
  elif i=8 and j=9 then return 0:
  elif i=8 and j=10 then return 0:
  elif i=8 and j=11 then return 0:
  elif i=8 and j=12 then return 0:
  elif i=8 and j=13 then return 0:
  elif i=8 and j=14 then return 0:
  elif i=8 and j=15 then return 0:
  elif i=8 and j=16 then return 2*x[8]:
  elif i=9 and j=10 then return 0:
  elif i=9 and j=11 then return 0:
  elif i=9 and j=12 then return 0:
  elif i=9 and j=13 then return 0:
  elif i=9 and j=14 then return -2*x[8]:
  elif i=9 and j=15 then return -x[6]:
  elif i=9 and j=16 then return -x[9]:
  elif i=10 and j=11 then return 0:
  elif i=10 and j=12 then return 0:
  elif i=10 and j=13 then return 0:
  elif i=10 and j=14 then return 0:
  elif i=10 and j=15 then return -2*x[8]:
  elif i=10 and j=16 then return 0:
  elif i=11 and j=12 then return 0:
  elif i=11 and j=13 then return 0:
  elif i=11 and j=14 then return 0:
  elif i=11 and j=15 then return 0:
  elif i=11 and j=16 then return x[11]:
  elif i=12 and j=13 then return 0:
  elif i=12 and j=14 then return -x[11]:
  elif i=12 and j=15 then return -x[10]:
  elif i=12 and j=16 then return -2*x[12]:
  elif i=13 and j=14 then return 0:
  elif i=13 and j=15 then return x[11]:
  elif i=13 and j=16 then return -x[13]:
  elif i=14 and j=15 then return 0:
  elif i=14 and j=16 then return 3*x[14]:
  elif i=15 and j=16 then return 2*x[15]:
  fi:
  if i=j then return 0 fi:  
  if j<i then return -L(j,i) fi:
  end:

  y := proc(i)
   if i=1 then return x[1]:
 elif i=2 then return 0*x[2]: 
 elif i=3 then return 1*x[3]:
 elif i=4 then return 2*x[4]:
 elif i=5 then return 2*x[5]:
 elif i=6 then return 3*x[6]:
 elif i=7 then return 1*x[7]:
 elif i=8 then return 4*x[8]:
 elif i=9 then return 3*x[9]:
 elif i=10 then return 4*x[10]:
 elif i=11 then return 5*x[11]:
 elif i=12 then return 4*x[12]:
 elif i=13 then return 5*x[13]:
 elif i=14 then return x[14]:
 elif i=15 then return 0*x[15]:
 elif i=16 then return 0*x[16]:
 fi: end: 

  E := diag(1,1,1,1,1,1,1,1,1,1,1,1,1,1,1,1):

  A := matrix(16,16):
  for i from 1 to 16 do
  for j from 1 to 16 do
  A[i,j] := L(j,i):
  od: od:

  f := proc(i)

  if i=1 then return matrix(1,16,[0,0,0,1,1,0,0,0,0,0,0,0,0,0,0,0]): 
  elif i=2 then return 
matrix(1,16,[0,0,0,0,0,x[13],0,-2*x[12],x[11],x[10],x[9],-2*x[8],x[6],
0,0,0]):
  elif i=3 then return 
matrix(1,16,[-x[11]*x[13]^2,0,x[13]*x[11]^2,x[13]^2*x[8]+
x[13]*x[11]*x[10]+x[12]*x[11]^2,0,0,-x[11]^3,-x[13]^2*x[5],
0,-x[11]*x[13]*x[5],-2*x[11]*x[12]*x[5]-3*x[11]^2*x[7]-
x[10]*x[13]*x[5]-x[13]^2*x[1]+2*x[13]*x[11]*x[3],-x[11]^2*x[5],
-x[11]*x[5]*x[10]+3*x[13]^2*x[14]+x[11]^2*x[3]-
2*x[13]*x[11]*x[1]-2*x[5]*x[13]*x[8],x[13]^3,0,0]):
  elif i=4 then return 
matrix(1,16,[-(x[13]*x[10]+2*x[12]*x[11])*x[13]^2,-x[13]^2*x[11]^2,
-(x[13]^2*x[8]-3*x[12]*x[11]^2-x[13]*x[11]*x[10])*x[13],
3*x[12]^2*x[11]^2+1/2*x[13]^3*x[6]+1/2*x[13]^2*x[9]*x[11]-
x[12]*x[13]^2*x[8]+x[13]^2*x[10]^2+3*x[13]*x[10]*x[12]*x[11],
0,0,3*x[11]*(x[13]^2*x[8]-x[12]*x[11]^2),-x[13]^2*(x[13]*x[3]-
3*x[11]*x[7]),0,1/2*x[13]*(-2*x[13]^2*x[1]+2*x[13]*x[11]*x[3]-
6*x[11]*x[12]*x[5]-3*x[10]*x[13]*x[5]),-3*x[10]*x[13]*x[5]*x[12]+
6*x[13]*x[11]*x[3]*x[12]+3*x[13]^2*x[7]*x[8]+x[10]*x[13]^2*x[3]-
9*x[11]^2*x[7]*x[12]+x[13]^3*x[16]-
2*x[13]^2*x[11]*x[2]-6*x[11]*x[12]^2*x[5]-2*x[13]^2*x[12]*x[1],
-3*x[13]^3*x[14]+x[13]^2*x[11]*x[1]+3*x[5]*x[13]^2*x[8]-
3*x[11]^2*x[5]*x[12],-3*x[11]*x[12]*x[10]*x[5]-
x[10]*x[13]^2*x[1]+6*x[13]^2*x[14]*x[12]+
2*x[13]^3*x[15]+3*x[11]^2*x[3]*x[12]-6*x[12]*x[13]*x[11]*x[1]-
6*x[13]*x[8]*x[5]*x[12]-x[3]*x[13]^2*x[8]+x[11]*x[13]^2*x[16],0,x[13]^4,
x[11]*x[13]^3]):
 elif i=5 then return 
simplify(multiply(matadd(f(4),multiply(f(3),-3*x[12]*E)),1/x[13]^2*E)):
  
    fi: end:

   k :=  proc(i)
    if i=1 then return f(1):
  elif i=2 then return f(2):
  elif i=3 then return 
simplify(multiply(matadd(f(4),multiply(f(3),-3*x[12]*E)),1/x[13]^2*E)):
  elif i=4 then return f(3):
  fi: end:

   r := proc(i)
     if i=1 then return 1/2*x[4]:
   elif i=2 then return 3*x[12]*x[8]-3/4*x[10]^2:
   elif i=3 then return -x[8]*x[13]^2-x[13]*x[11]*x[10]-x[12]*x[11]^2:
  fi: end:

   M := proc(i)
if i=1 then return 
matadd(matadd(k(3),multiply(k(2),r(1)*E)),multiply(k(1),r(2)*E)):
elif i=2 then return matadd(k(4),multiply(k(1),r(3)*E)):
fi: end:

    p := proc(i)
 
    if i=1 then return add(k(1)[1,j]*x[j],j=1..16):
   elif i=2 then return add(k(2)[1,j]*x[j],j=1..16):
   elif i=3 then return add(M(1)[1,j]*x[j],j=1..16):
   elif i=4 then return add(M(2)[1,j]*x[j],j=1..16):
    fi: end:

     P := proc(i,j::posint)
if i=1 and j<17 then return simplify(add(diff(p(1),x[l])*L(j,l),l=1..16)):
elif i=2 and j<17 then return simplify(add(diff(p(2),x[l])*L(j,l),l=1..16)):
elif i=3 and j<17 then return simplify(add(diff(p(3),x[l])*L(j,l),l=1..16)):
elif i=4 and j<17 then return simplify(add(diff(p(4),x[l])*L(j,l),l=1..16)):
  fi: end:
\end{verbatim}

We add the Gap4 program to compute the bracket L(i,j) which is equal to Bg[i]*Bg[j] in 
the Gap4 program.

\begin{verbatim}
 L := SimpleLieAlgebra("F",4,Rationals);;
 R := RootSystem(L);;
 P := PositiveRoots(R);;
 x := PositiveRootVectors(R);;
 y := NegativeRootVectors(R);;
 e := x[11]+x[12]+x[13];;
 IsNilpotentElement(L,e);;
 if true then FindSl2(L,e);;
 Bs := BasisVectors(Basis(FindSl2(L,e)));;
 F := function(i)
 return (i)*((Bs[1]*Bs[3])*Bs[1]);;
 end;;

 numbers := [1..20];;
 for i in numbers do
   if F(i) = (2)*e then
     f := (i)*Bs[3];;
     fi;;
     od;;
     
 h := e*f;;

 g := LieCentralizer(L,Subspace(L,[e]));;
 Bg := BasisVectors(Basis(g));;
 z := LieCentre(g);;
 Bz := BasisVectors(Basis(z));;
 fi;;
\end{verbatim}

\backmatter

\bibliographystyle{smfplain}

\end{document}